\documentclass[aos,MSNbibl,dvips]{arximspdf}
\usepackage{graphicx}
\doi{10.1214/13-AOS1174} 
\volume{41}
\issue{6}
\pubyear{2013}
\firstpage{3050}
\lastpage{3073}

\makeatletter
\newcommand{\rrVert}{\Vert}
\newcommand{\llVert}{\Vert}
\newcommand{\E}{\mathrm{E}}
\newcommand{\cov}{\operatorname{Cov}}

\newcommand{\bZ}{\mathbb{Z}}
\newcommand{\bR}{\mathbb{R}}

\newcommand{\X}{X_t}

\newcommand{\BEL}{\mathrm{BEL}}

\newtheorem{theorem}{Theorem}

\newtheorem{Corollary}{Corollary}
\newtheorem{Lemma}{Lemma}
\newproclaim{Assumption}{Assumption}
\newproclaim{Fact}{Fact}

\newproclaim{Definition}{Definition}
\newproclaim{ass}{Assumptions}
\newproclaim{remark}{Remark}
\makeatother

\begin{document}
\begin{frontmatter}

\title{A nonstandard empirical likelihood for time series}
\runtitle{Nonstandard EL}

\begin{aug}
\author[a]{\fnms{Daniel J.} \snm{Nordman}\corref{}\thanksref{t1}\ead[label=e1]{dnordman@iastate.edu}},
\author[b]{\fnms{Helle} \snm{Bunzel}\ead[label=e2]{hbunzel@iastate.edu}}
\and\break
\author[c]{\fnms{Soumendra N.} \snm{Lahiri}\thanksref{t2}\ead[label=e3]{snlahiri@ncsu.edu}}
\thankstext{t1}{Supported in part by NSF DMS-09-06588.}
\thankstext{t2}{Supported in part by NSF DMS-10-07703 and NSA
H98230-11-1-0130.}
\runauthor{D.~J. Nordman, H. Bunzel and S.~N. Lahiri}
\affiliation{Iowa State University, Iowa State University and Aarhus University,
and North
Carolina State University}
\address[a]{D.~J. Nordman\\
Department of Statistics\\
Iowa State University\\
Ames, Iowa 50011\\
USA\\
\printead{e1}}
\address[b]{H. Bunzel\\
Department of Economics\\
Iowa State University\\
Ames, Iowa 50011\\
USA\\
and\\
CREATES\\
Aarhus University\\
Aarhus, DK-8000\\
Denmark\\
\printead{e2}}

\address[c]{S.~N. Lahiri\\
Department of Statistics\\
North Carolina State University\\
Raleigh, North Carolina 27695-8203\\
USA\\
\printead{e3}}
\end{aug}

\received{\smonth{4} \syear{2012}}
\revised{\smonth{7} \syear{2013}}

%
\begin{abstract}
Standard blockwise empirical likelihood (BEL) for stationary, weakly
dependent time series
requires specifying a fixed block length as a tuning parameter for
setting confidence regions.
This aspect can be difficult and impacts coverage accuracy.
As an alternative,
this paper proposes a new version of BEL based on a simple, though
nonstandard, data-blocking rule which uses a data block of every
possible length.
Consequently, the method does not involve the usual block selection
issues and is also anticipated to exhibit
better coverage performance. Its nonstandard blocking scheme, however,
induces nonstandard asymptotics and requires a significantly
different development
compared to standard BEL.
We establish the large-sample distribution of log-ratio statistics
from the new BEL method
for calibrating confidence regions for mean or
smooth function parameters of time series.
This limit law is not the usual chi-square one,
but is distribution-free and
can be reproduced through straightforward simulations. Numerical
studies indicate that the proposed method generally exhibits better
coverage accuracy
than standard BEL.
\end{abstract}

%
\pdfkeywords{62G09, 62G20, 62M10, Brownian motion, confidence regions, stationarity, weak dependence}
\end{frontmatter}

\section{Introduction}\label{sec1}
For independent, identically distributed data (i.i.d.), Owen \cite
{O88,O90} introduced empirical likelihood (EL)
as a general methodology for re-creating
likelihood-type inference
without a joint distribution for the data, as typically specified in
parametric likelihood.
However, the i.i.d. formulation of EL fails for dependent data by
ignoring the underlying dependence structure. As a remedy,
Kitamura \cite{K97} proposed so-called blockwise empirical likelihood
(BEL) methodology
for stationary, weakly dependent processes, which has been shown to
provide valid inference in various scenarios with time series
(cf. \cite{B05,B09,CW09,LZ11,NL07,WC11}). Similarly to the i.i.d. EL
version, BEL creates
an EL log-ratio statistic having a chi-square limit for inference, but
the BEL construction crucially involves blocks of consecutive
observations in time, rather than individual observations. This
data-blocking serves to capture the underlying time dependence and
related concepts have also proven important in defining resampling
methodologies for dependent data,
such as block bootstrap \cite{H85,K89,LS92} and time subsampling
methods \cite{C86,PR93,PRW99}. However, the coverage accuracy
of BEL can depend crucially on the block length selection, which is a
fixed value $1 \leq b \leq n$
for a given sample size~$n$, and appropriate choices can vary with the
underlying process (a point
briefly illustrated
at the end of this section).

To advance the BEL methodology in a direction away from block selection
with a goal of
improved coverage accuracy, we propose an alternative version of BEL
for stationary, weakly dependent time series, called
an expansive block empirical likelihood (EBEL). The EBEL method
involves a nonstandard, but simple, data-blocking rule where a data block
of every possible length is used. Consequently, the method does not
involve a block length choice in the standard sense. We investigate
EBEL in the prototypical problem
of inference about the process mean or a smooth function of means. For
setting confidence regions for such parameters, we establish the limiting
distribution of log-likelihood ratio statistics from the EBEL method.
Because of the nonstandard blocking scheme,
the justification of this limit distribution
requires a new and substantially different treatment compared
to that of standard BEL (which closely resembles that of EL for i.i.d. data
in its large-sample development; cf. \cite{O90,QL94}).
In fact, unlike with standard BEL or EL for i.i.d. data, the limiting
distribution involved is nonstandard and \textit{not} chi-square.
However, the EBEL limit law is distribution-free,
corresponding to a special integral of standard Brownian motion on
$[0,1]$, and so can be easily
approximated through simulation to obtain appropriate quantiles
for calibrating confidence regions.
In addition, we anticipate that the EBEL method
may have generally better coverage accuracy than standard BEL methods,
though formally
establishing and comparing convergence rates is beyond the scope of
this manuscript (and, in fact,
optimal rates and block sizes for
even standard BEL remain to be determined).
Simulation studies, though, suggest that
interval estimates from the EBEL method can perform much better than
the standard BEL approach, especially
when the later employs a poor block choice, and be less sensitive to
the dependence strength of the underlying process.

The rest of manuscript is organized as follows. We end this section by
briefly recalling
the standard BEL construction with overlapping blocks and its
distributional features.
In Section~\ref{sec2}, we separately describe the EBEL method for inference on
process means and smooth function model parameters, and
establish the main distributional results in both cases. These results
require introducing a new type of limit
law based on Brownian motion, which is also given in Section~\ref{sec2}.
Additionally, Section~\ref{sec2.1} describes
how the usual EL theory developed by Owen \cite{O88,O90}, and often
underlying many EL arguments including
the time series extensions of BEL \cite{K97}, fails here and requires
new technical developments;
consequently, the theory provided may be useful for future
developments of EL (with an example given in Section~\ref{sec2.4}).
Section~\ref{sec3} provides a numerical study of
the coverage accuracy of the EBEL method and comparisons to standard
BEL. Section~\ref{sec4} offers some concluding remarks and heuristic arguments
on the expected performance
of EBEL. Proofs of the main results appear in Section~\ref{sec5} and in
supplementary materials \cite{NBL13}, where the latter also presents
some additional simulation summaries.

\begin{figure}

\includegraphics{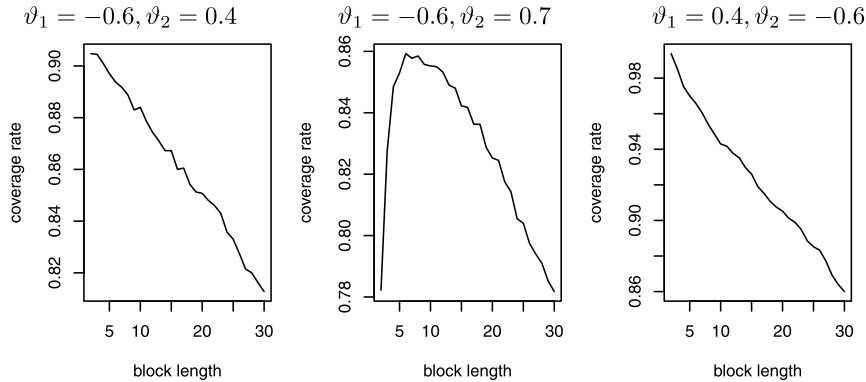}

\caption{Plot of coverage rates for 90\% BEL intervals for the process
mean $\E X_t=\mu$ over various
blocks $b=2,\ldots,30$, based on samples of size $n=100$ from three
$\operatorname{MA}(2)$ processes
$X_t = Z_t + \vartheta_1 Z_{t-2}+\vartheta_2$ with i.i.d. standard
normal innovations $\{Z_t\}$ (from $4000$ simulations).}\label{fig1}
\end{figure}

To motivate what follows, we briefly recall the BEL construction,
considering, for concreteness,
inference about the mean $\E X_t=\mu\in\mathbb{R}^d$ of a
vector-valued stationary stretch
$X_1,\ldots,X_n$.
Upon choosing an integer block length $1 \leq b \leq n$, a collection
of maximally overlapping (OL) blocks of length $b$ is given by $\{
(X_i,\ldots,X_{i+b-1})\dvtx i=1,\ldots,N_b \equiv n-b+1\}$. For a
given $\mu
\in\mathbb{R}^d$ value, each block
in the collection provides a centered block sum $B_{i,\mu} \equiv\sum
_{j=i}^{i+b-1} (X_j-\mu)$ for defining a BEL function
%
%
\begin{equation}
\label{eqn:stan} L_{\BEL,n}(\mu) = \sup\Biggl\{ \prod
_{i=1}^{N_b} p_i \dvtx p_i
\geq0, \sum_{i=1}^{N_b}p_i = 1,
\sum_{i=1}^{N_b}p_i
B_{i,\mu} = 0_d \Biggr\}
\end{equation}
and corresponding BEL ratio $R_{\BEL,n}(\mu) = L_n(\mu)/N_b^{-N_b}$,
where above $0_d=(0,\ldots,0)^\prime\in\mathbb{R}^d$. The function
$L_{\BEL,n}(\mu)$ assesses the plausibility of a value $\mu$ by
maximizing a multinomial
likelihood from probabilities $\{p_i\}_{i=1}^{N_b}$
assigned to the centered block sums $B_{i,\mu}$
under a zero-expectation constraint.
Without the linear mean constraint in (\ref{eqn:stan}), the multinomial
product is maximized when each $p_i=1/N_b$ (i.e., the empirical
distribution on blocks), defining the ratio
$R_{\BEL,n}(\mu)$. Under certain mixing and moment conditions entailing
weak dependence, and if the block $b$
grows with the sample size $n$ but at a smaller rate (i.e.,
$b^{-1}+b^2/n\rightarrow0$ as $n\rightarrow\infty$), the log-EL ratio
of the standard BEL
has chi-square limit
%
%
\begin{equation}
\label{eqn:stanl} -\frac{2}{b} \log R_{\BEL,n}(\mu_0)
\stackrel{d} {\rightarrow} \chi^2_d,
\end{equation}
at the true mean parameter $\mu_0\in\mathbb{R}^d$; cf. Kitamura
\cite
{K97}. Here $b^{-1}$ represents
an adjustment in (\ref{eqn:stanl}) to account for OL blocks and, for
i.i.d. data, a block length $b=1$ above produces the EL distributional
result of Owen \cite{O88,O90}. To illustrate the connection between
block selection
and performance, Figure~\ref{fig1} shows the coverage rate of nominal 90\% BEL
confidence intervals
$\{\mu\in\mathbb{R}\dvtx-2/b \log R_{\BEL,n}(\mu_0) \leq\chi
^2_{1,0.9}\}$,
as a function of the
block size $b$,
for estimating the mean of three different $\operatorname{MA}(2)$ processes
based on samples of size $n=100$. One observes that the coverage
accuracy of BEL varies
with the block length and that the best block size can depend on
the underlying process.
The EBEL method described next is a type of modification of the OL BEL
version, without a particular fixed
block length selection $b$.

%

%
%
%
%

\section{Expansive block empirical likelihood}\label{sec2}

\subsection{Mean inference}\label{sec2.1}

Suppose $X_1,\ldots,X_n$ represents a sample from a strictly
stationary process $\{\X\dvtx t\in\bZ\}$ taking values in $\bR^d$, and
consider a problem about inference
on the process mean $\E X_t=\mu\in\mathbb{R}^d$. While the BEL uses
data blocks of a fixed length $b$ for a given sample size $n$,
the EBEL uses overlapping data blocks $\{ (X_1), (X_1,X_2),\ldots,
(X_1,\ldots,X_n)\}$ that vary in length up to the longest block
consisting of
the entire time series. Hence this block collection, which constitutes
a type of forward ``scan'' in the block subsampling
language of McElroy and Politis~\cite{MP07}, contains a data block of
every possible length $b$
for a given sample size $n$. This block sequence also appears in
fixed-$b$ asymptotic schemes
\cite{KV02} and related self-normalization approaches; cf. Shao \cite
{S10}, Section~2; see also Section~\ref{sec4} here. In this sense, these blocks
are interesting and novel to consider in a BEL framework. Other block
schemes may be possible
and potentially applied for practical gain (e.g., improved power), where the
theoretical results of this paper could also directly apply.
We leave this largely for future research, but we shall give one
example of a modified, though related, blocking scheme in Section~\ref{sec2.4}
while focusing the exposition on the block collection above [i.e., the
alternative blocking incorporates
a backward scan $\{ (X_n), (X_n,X_{n-1}),\ldots, (X_n,\ldots,X_1)\}$
with similar theoretical development].

Let $w\dvtx[0,1]\rightarrow[0,\infty)$ denote a nonnegative
weighting function.
To assess the likelihood of a given value of $\mu$, we create centered
block sums $T_{i,\mu} = w(i/n) \sum_{j=1}^{i}(X_j - \mu)$,
$i=1,\ldots
,n$, and define a EBEL function
%
%
\begin{equation}
\label{eqn:EBEL} L_{n}(\mu) = \sup\Biggl\{ \prod
_{i=1}^{n} p_i \dvtx p_i
\geq0, \sum_{i=1}^{n}p_i = 1,
\sum_{i=1}^{n} p_i
T_{i,\mu} = 0_d \Biggr\}
\end{equation}
and ratio $R_n(\mu)= n^{-n}L_n(\mu) $.
After defining the block sums, the computation of $L_{n}(\mu)$ is
analogous to the BEL version and essentially the same as that described
by Owen \cite{O88,O90} for i.i.d. data. Namely, when the zero $0_d$
vector lies in the interior convex hull of $\{T_{i,\mu}\dvtx
i=1,\ldots,n\}
$, then $L_n(\mu)$ is the uniquely achieved maximum at probabilities
$p_i= 1/[n(1 + \lambda_{n,\mu}^\prime T_{i,\mu})]>0$, $i=1,\ldots,n$,
with a Lagrange multiplier $\lambda_{n,\mu}\in\mathbb{R}^d$ satisfying
%
%
\begin{equation}
\label{eqn:t0} \sum_{i=1}^n
\frac{T_{i,\mu}}{n(1 + \lambda_{n,\mu}^\prime T_{i,\mu})
} = 0_d;
\end{equation}
see \cite{O90} for these and other computational details. Regarding the
weight function above in the EBEL formulation, more details are
provided below and in Section~\ref{sec2.2}.

The next section establishes the limiting distribution of the log-EL
ratio from the EBEL method for setting confidence regions for the
process mean $\mu$ parameter. However, it is helpful to initially describe
how the subsequent developments of EL differ from previous ones with
i.i.d. or weakly dependent data (cf. \cite{K97} for BEL).
The standard arguments for developing EL results, due to Owen \cite
{O90} (page 101), typically begin from algebraically re-writing (\ref{eqn:t0})
to expand the Lagrange multiplier. If we consider the real-valued case
$d=1$ for simplicity, this becomes
\[
\lambda_{n,\mu} = \frac{\sum_{i=1}^n T_{i,\mu}}{\sum_{i=1}^n
T_{i,\mu
}^2 } + \frac{ \lambda_{n,\mu}^2}{\sum_{i=1}^n T_{i,\mu}^2} \sum
_{i=1}^n \frac{T_{i,\mu}^3}{1 + \lambda_{n,\mu}^\prime T_{i,\mu}}.
\]
In the usual independence or weak dependence cases of EL [e.g., where
$B_{i,\mu}$ from~(\ref{eqn:stan}) replaces $T_{i,\mu}$ in the Lagrange
multiplier above], the first right-hand side term dominates the second,
which gives
a substantive form for $\lambda_{n,\mu}$ as a ratio of sample means
and consequently drives the large sample results (i.e., producing
chi-square limits). However, in the EBEL
case here, both terms on the right-hand side above have the \textit{same}
order, implying that the standard approach to developing EL results
breaks down under the EBEL blocking scheme. The proofs here use a
different EL argument than the standard one mentioned above \cite{K97,O90}, involving no asymptotic expansions of the Lagrange multiplier or
Taylor expansions of
the EL ratio based on these.

The large sample results for the EBEL method require two mild
assumptions stated below. Let $\mathcal{C}_d[0,1]$
denote the metric space of all $\mathbb{R}^d$-valued continuous
functions on $[0,1]$ with the supremum metric $\rho(g_1,g_2)\equiv
\sup_{0 \leq t \leq1}\|g_1(t)-g_2(t)\|$, and let $B(t) =
(B_1(t),\ldots
,B_d(t))^\prime$, $0 \leq t \leq1$, denote
a $\mathcal{C}_d[0,1]$-valued random variable where $B_1(t),\ldots
,B_d(t)$ are i.i.d. copies of standard Brownian motion on $[0,1]$.

\begin{ass*}
\begin{longlist}[((A.2))]
\item[(A.1)] The weight function $w\dvtx[0,1]\rightarrow[0,\infty)$ is
continuous on $[0,1]$ and is strictly positive on an interval $(0,c)$
for some $c \in(0,1]$.
\item[(A.2)] Let $\E X_t =\mu_0 \in\mathbb{R}^d$ denote the true mean
of the stationary process $\{X_t\}$ and suppose $d\times d$ matrix
$\Sigma= \sum_{j=-\infty}^\infty\cov(X_0,X_j)$ is positive definite.
For the empirical process $S_n(t)$ on $t\in[0,1]$ defined by linear
interpolation of $\{ S_n(i/n) = \sum_{j=1}^i (X_j - \mu_0) \dvtx
i=0,\ldots
, n \}$ with $S_n(0)=0$, it holds that
$ S_n(\cdot)/n^{1/2} \stackrel{d}{\rightarrow} \Sigma^{1/2} B(\cdot)$
in $\mathcal{C}_d[0,1]$.
\end{longlist}
\end{ass*}

Assumption (A.1) is used to guarantee that, in probability, the EBEL
ratio $R_n(\mu_0)$
positively exists at the true mean, which holds for uniformly weighted
blocks $w(t)=1$, $t\in[0,1]$, for example. Assumption (A.2) is a
functional central limit theorem for weakly dependent data, which
holds under appropriate mixing and moment conditions on $\{X_t\}$ \cite{H84}.

\subsection{Main distributional results}\label{sec2.2}

To state the limit law for the log-EBEL ratio~(\ref{eqn:EBEL}), we
require a result
regarding a vector $B(t) = (B_1(t),\ldots,B_d(t))^\prime$, $0 \leq t
\leq1$, of i.i.d. copies $B_1(t),\ldots,B_d(t)$ of standard Brownian
motion on $[0,1]$. Indeed, the limit distribution
of $-2 \log R_n(\mu_0)$ is a nonstandard functional of the vector of
Brownian motion $B(\cdot)$.
Theorem \ref{theorem1} identifies key elements of the limit law and describes some
of its basic structural properties.

%
\begin{theorem} \label{theorem1}
Suppose that $B(t) = (B_1(t),\ldots,B_d(t))^\prime$, $0 \leq t \leq
1$, is defined on a probability space,
and let $f(t)=w(t)B(t)$, $0\leq t\leq1$, where $w(\cdot)$ satisfies
assumption \textup{(A.1)}. Then, with probability 1 (w.p.1), there
exists an $\mathbb{R}^d$-valued random vector $Y_d$ satisfying the following:
\begin{longlist}[(iii)]
\item[(i)] $Y_d$ is the unique minimizer of
\[
g_d(a) \equiv-\int_0^1 \log
\bigl( 1 + a^\prime f(t)\bigr)\,dt\qquad \mbox{for $a\in\overline{K}_d$},
\]
where $\overline{K}_d \equiv\{ y \in\mathbb{R}^d \dvtx\min_{0
\leq t
\leq1} (1 + y^\prime f(t)) \geq0 \}$ is the closure of
$K_d \equiv\{ y \in\mathbb{R}^d \dvtx\min_{0 \leq t \leq1} (1 +
y^\prime f(t)) > 0 \}$; the latter set
is open, bounded and convex in $\mathbb{R}^d$ (w.p.1). On $K_d$, $g_d$
is also real-valued, strictly convex and infinitely differentiable (w.p.1).
\item[(ii)] $-\infty< g_d(Y_d) <0$, $ {Y_d^\prime\int_0^1 f(t)\,dt >0}$,
$ {0 \leq\int_0^1 \frac{Y_d^\prime f(t)}{ 1 + Y_d^\prime f(t)}\,dt <
\infty}$.
\item[(iii)] If $Y_d \in K_d$, then $Y_d$ is the unique solution to $
{\int_0^1 \frac{f(t)}{1 + a^\prime f(t)}\,dt = 0_d}$ for $a\in K_d$, and
if $ {\int_0^1 \frac{f(t)}{1 + a^\prime f(t)}\,dt = 0_d}$ has a solution
$a\in K_d$, then this solution is uniquely $Y_d$.
\end{longlist}
\end{theorem}
We use the subscript $d$ in Theorem \ref{theorem1} to denote the dimension
of either the random vector $Y_d$, the space $K_d$ or the arguments of
$g_d$. The function $g_d$ is well defined and convex on $\overline
{K}_d$, though possibly
$g_d(a) = +\infty$ for some $a\in\partial K_d = \{ y \in\mathbb{R}^d
\dvtx\min_{0 \leq t \leq1} (1 + y^\prime f(t)) = 0 \}$ on the
boundary of
$K_d$; a minimizer of $g_d(\cdot)$ may also occur on $\partial K_d
\cap
\{y \in\mathbb{R}^d\dvtx g_d(y)\leq0 \}$. Importantly, the probability
law of $g_d(Y_1)$ is distribution-free, and because standard Brownian motion
is fast and straightforward to simulate, the distribution of $g_d(Y_d)$
can be approximately numerically. Parts (ii) and (iii) provide
properties for characterizing and identifying the minimizer $Y_d$. For example,
considering the real-valued case $d=1$, it holds that $K_1 = (m, M)$
where $m= -[\max_{0 \leq t \leq1} f(t)]^{-1} < 0 < M = - [\min_{0
\leq
t \leq1} f(t)]^{-1}$ and the derivative $d g_1(a)/da$ is strictly
increasing on $K_1$ by convexity. Because the derivative of $g_1$ at 0
is $- \int_0^1 f(x) \,dx$, parts (ii)--(iii) imply that if $- \int_0^1
f(x) \,dx<0$, then either $Y_1=m$ or $Y_1$ solves
$d g_1(a)/da =0$ on $m < a \leq0$; alternatively, if $- \int_0^1 f(x)
\,dx>0$, then $Y_1=M$ or $Y_1$ solves
$d g_1(a)/da =0$ on $0 \leq a < M$. Additionally, while the weight
function $w(\cdot)$ influences
the distribution of $g_d(Y_d)$, the scale of $w(\cdot)$ does not;
defining $f$ with $w$ or $c w$, for a nonzero $c\in\mathbb{R}$,
produces the same
minimized value $g_d(Y_d)$.

We may now state the main result on the large-sample behavior of the
EBEL log-ratio evaluated at the true process mean $\E X_t =\mu_0 \in
\mathbb{R}^d$. Recall that, when $L_n(\mu_0)>0$ in (\ref{eqn:EBEL}),
the EBEL log-ratio admits a representation (\ref{eqn:t0})
at $\mu_0$
in terms of the Lagrange multiplier $\lambda_{n,\mu_0} \in\mathbb{R}^d$.

\begin{theorem} \label{theorem2}
Under assumptions \textup{(A.1)--(A.2)}, as $n \rightarrow\infty$:
\begin{longlist}[(ii)]
\item[(i)] $ n^{1/2} \Sigma^{1/2} \lambda_{n,\mu_0} \stackrel
{d}{\rightarrow} Y_d $;
\item[(ii)] $ {-\frac{1}{n} \log R_n(\mu_0) \stackrel{d}{\rightarrow}
-g_d(Y_d)}$,
\end{longlist}
for $Y_d$ and $g_d(Y_d)$ defined as in Theorem \ref{theorem1}, and $\Sigma= \sum
_{j=-\infty}^\infty\cov(X_0,X_j)$.
\end{theorem}
From Theorem \ref{theorem2}(i), the Lagrange multiplier in the EBEL method
has a limiting distribution which is not the typical normal one, as in
the standard BEL case. This has a direct impact on the limit law of the
EBEL ratio statistic.
As Theorem \ref{theorem2}(ii) shows, the negative logarithm of the EBEL ratio
statistic, scaled\vadjust{\goodbreak} by the inverse of the sample size, has a nonstandard
limit, given by the functional $-g_d(Y_d)$ of the vector of Brownian motion
$B(\cdot)$ (cf.~Theorem \ref{theorem1}), that critically depends on the limit $Y_d$
of the scaled Lagrange multiplier.
The distribution of $-g_d(Y_d)$ is free of any population parameters so
that quantiles of $-g_d(Y_d)$, which are easy to compute numerically,
can be used to
calibrate the EBEL confidence regions.
As $-g_d(Y_d)$ is a strictly positive random variable, an approximate
$100(1-\alpha)\%$ confidence region for $\mu_0$
can be computed as
\[
\bigl\{\mu\in\mathbb{R}^d \dvtx- n^{-1} \log
R_n(\mu_0) \leq a_{d,1-\alpha}\bigr\},
\]
where $a_{d,1-\alpha}$
is the lower $(1-\alpha)$ percentile of $-g_d(Y_d)$. When $d=1$, the
confidence region is an interval; for $d>2$, the region is
guaranteed to be connected without voids in $\mathbb{R}^d$. In
contrast to the standard BEL (\ref{eqn:stanl}), EBEL confidence regions
do not require a similar fixed choice of block size.

We next provide additional results that give the limit distribution of
the log-EBEL ratio statistic under a sequence of
local alternatives and that also show the size of a EBEL confidence
region will be no larger than
$O_p(n^{-1/2})$ in diameter around the true mean $\E X_t = \mu_0$. Let
%
%
\begin{equation}
\label{eqn:Gn} G_n \equiv\bigl\{ \mu\in\mathbb{R}^d \dvtx
R_n(\mu)\geq R_n(\mu_0)>0\bigr\}
\end{equation}
be the collection
of mean parameter values which are at least as likely as $\mu_0$, and
therefore elements of a EBEL confidence region whenever the true mean is.
%
%
\begin{Corollary}\label{co1} Suppose the assumptions of Theorem \ref{theorem2} hold. For
$c \in\mathbb{R}^d$, define $f_{c}(t) =w(t)[B(t) + t\Sigma^{-1/2}
c]$, $t\in[0,1]$, in terms of the vector of Brownian motion $B(t)$.
\begin{longlist}[(ii)]
\item[(i)] Then, as $n\rightarrow\infty$, $-n^{-1}\log R_n(\mu_0 \pm
n^{-1/2}c) \stackrel{d}{\rightarrow}$
\[
-\min\biggl\{ -\int_0^1 \log\bigl(1 +
a^\prime f_{c}(t)\bigr)\,dt\dvtx a \in\mathbb{R}^d,
\min_{0 \leq t \leq1} \bigl(1 + a^\prime f_{c}(t)\bigr)
\geq0 \biggr\};
\]
\item[(ii)] $ \sup\{\|\mu-\mu_0\|\dvtx\mu\in G_n\} =O_p(n^{-1/2})$, for
$G_n$ in (\ref{eqn:Gn}).
\end{longlist}
\end{Corollary}

Hence along a sequence of local alternatives ($n^{-1/2}$ away from the
true mean), the log-EBEL ratio converges to a random variable, defined
as the optimizer of an integral involving Brownian motion; this
resembles Theorem \ref{theorem1} [involving $f(t)=w(t)B(t)$ there], but the
integrated function $f_c(\cdot)$ has an addition term $w(t) t \Sigma
^{-1/2} c$ under the alternative. With respect to Corollary \ref{co1}(i), the
involved limit distribution can be described with similar properties as
in Theorem \ref{theorem1} upon replacing $f(t)$ with $f_c(t)$ there.
In particular, the limiting distribution
under the scaled alternatives depends on $\Sigma^{-1/2} c$, similarly
to the normal theory case
(e.g., with standard BEL) where $\Sigma^{-1/2} c$ determines the noncentrality
parameter of a noncentral chi-square distribution.\vadjust{\goodbreak}

We note that Theorem \ref{theorem2} remains valid for potentially
negative-valued weight functions $w(\cdot)$ as well. Simulations have
shown that, with weight functions oscillating between positive and
negative values on $[0,1]$ [e.g., $w(t)=\sin(2\pi t)$],
EBEL intervals for the process mean perform consistently well in terms
of coverage accuracy. However, with weight functions $w(\cdot)$ that
vary in sign, a result as in Corollary \ref{co1}(ii) fails to hold. Hence, the
weight functions $w(\cdot)$ considered are nonnegative
as stated in assumption (A.1).

\begin{remark}\label{re1}
The EBEL results in Theorem \ref{theorem2}
also extend to certain parameters described by general estimating
functions; for examples and similar EL results in the i.i.d. and time
series cases, respectively, see \cite{QL94} and \cite{K97}. Suppose
$\theta\in\mathbb{R}^p$ represents a parameter of interest and
$G(\cdot;\cdot) \in\mathbb{R}^d \times\mathbb{R}^p \rightarrow
\mathbb
{R}^{p} $ is a vector of $p$ estimating functions such that $\E G(X_t;
\theta_0)=0_p$ holds at the true parameter value $\theta_0$.
The previous process mean case corresponds to $G(X_t; \mu) = X_t-\mu$
with $X_t, \mu\in\mathbb{R}^d$, $d=p$. A EBEL ratio statistic
$R_n(\theta)$ for $\theta$ results by replacing $T_{u,i} =w(i/n) \sum
_{j=1}^i (X_j-\mu)$ and $0_d$ with $T_{\theta,i} =w(i/n) \sum_{j=1}^i
G(X_j;\theta)$ and $0_p$ in (\ref{eqn:EBEL}). Under the conditions of
Theorem \ref{theorem2} [substituting $G(X_j;\theta_0)$ for $X_j-\mu_0$
in assumption (A.2)],
\[
-\frac{1}{n} \log R_n(\theta_0) \stackrel{d} {
\rightarrow} -g_p(Y_p)
\]
holds as $n\rightarrow\infty$ with $Y_p$ and $g_p(Y_p)$ as defined in
Theorem \ref{theorem1}, generalizing
Theorem \ref{theorem2} and following by the same proof. The next
section considers extensions of the EBEL approach to a different class
of time series parameters.
\end{remark}

\subsection{Smooth function model parameters}\label{sec2.3}
We next consider extending the EBEL method for inference on a broad
class of parameters
under the so-called ``smooth function model;'' cf. \cite{BG78,H92}.
For independent and time series data, respectively,
Hall and La
Scala \cite{HL90} and Kitamura \cite{K97} have considered EL inference
for similar parameters; see also \cite{O90}, Section~4.

If $\E X_t = \mu_0 \in\mathbb{R}^d$ again denotes the
true mean of the process, the target parameter of interest is given by
%
%
\begin{equation}
\label{eqn:H} \theta_0 = H(\mu_0) \in
\mathbb{R}^p,
\end{equation}
based on a smooth function $H(\mu) = (H_1(\mu),\ldots,H_p(\mu
))^\prime$
of the mean parameter $\mu$, where $H_i\dvtx\mathbb{R}^d\rightarrow
\mathbb
{R}$ for $i=1,\ldots,p$ and $p \leq d$. This framework allows a large
variety of parameters
to be considered such as sums, differences, products and ratios of
means, which can be used, for example,
to formulate parameters such as covariances and autocorrelations
as functions of the $m$-dimensional moment structure (for a fixed $m$)
of a time series.
For a univariate stationary series $U_1,\ldots,U_n$, for instance, one
can define a multivariate
series $X_t$ based on transformations of $(U_t,\ldots,U_{t+m-1})$ and
estimate parameters for
the process $\{U_t\}$ based on appropriate functions $H$ of the mean of $X_t$.
The correlations $\theta_0 = H(\mu_0)$ of $\{U_t\}$ at lags $m$ and
$m_1<m$, for example, can be formulated
in (\ref{eqn:H}) by $H(x_1,x_2,x_3,x_4) = (x_3-x_1^2,
x_4-x_1^2)^\prime
/[x_2-x_1^2]$ and $\E X_t = \mu_0$ for $X_t=(U_t,U_t^2, U_{t}U_{t+m_1},
U_{t}U_{t+m})^\prime\in\mathbb{R}^4$.
\cite{K89} and \cite{L03} (Chapter~4) provide further examples of
smooth function parameters.

For inference on the parameter $\theta= H(\mu)$, the EBEL ratio is
defined as
\[
R_n(\theta) = \sup\Biggl\{ \prod_{i=1}^{n}
p_i \dvtx p_i \geq0, \sum_{i=1}^{n}p_i
= 1, \sum_{i=1}^{n} p_i
T_{i,\mu} = 0_d, \mu\in\mathbb{R}^d, H(\mu)=
\theta\Biggr\},
\]
and its limit distribution is provided next.

\begin{theorem} \label{theorem13}
In addition to the assumptions of Theorem \ref{theorem2}, suppose $H$ from~(\ref
{eqn:H}) is continuously differentiable in a neighborhood of $\mu_0$
and that $\nabla_{\mu_0}$ has rank $p \leq d$, where $\nabla_\mu
\equiv
[\partial H_i(\mu)/\partial\mu_j]_{i=1,\ldots,p; j=1,\ldots,d}$
denotes the $p \times d$ matrix
of first-order partial derivatives of $H$. Then, at the true parameter
$\theta_0=H(\mu_0)$, as $n\rightarrow\infty$,
\[
-\frac{1}{n} \log R_{n}(\theta_0) \stackrel{d} {
\rightarrow} -g_p(Y_p)
\]
with $Y_p$ and $g_p(Y_p)$ as defined in Theorem \ref{theorem1}.
\end{theorem}

Theorem \ref{theorem13} shows that the log-EBEL ratio statistic for the parameter
$\theta_0=H(\mu_0)\in\mathbb{R}^p$ under
the smooth function model continues to have a limit of the same form as
that in the case
of the EBEL for the mean parameter $\mu_0 \in\mathbb{R}^d$ itself. The
main difference
is that the functional $g_p(Y_p)$ is now defined in terms of a
$p$-dimensional Brownian motion
as in Theorem \ref{theorem1}, but with $p \leq d$, where $p$ denotes the dimension
of the parameter $\theta_0$; see also Remark \ref{re1}. It is interesting to
note that, similarly to the traditional profile likelihood theory in a
parametric set-up with i.i.d. observations, the limit law here does not
depend on the function $H$ as long as the matrix $\nabla_{\mu_0}$ of
the first order partial derivatives of $H$ at $\mu=\mu_0$ has full
rank~$p$.
Due to the nonstandard blocking, the proof of this EBEL result again
requires a different
development compared to the one for standard BEL (cf. \cite{K97}) that mimics
the i.i.d. EL case (cf. \cite{HL90,O90}) involving expansion of
Lagrange multipliers.

\subsection{Extensions to other data blocking}\label{sec2.4}
As mentioned in Section~\ref{sec2.1}, other versions of EBEL may be possible
with other data blocking schemes,
which likewise involve no fixed block selection in the usual BEL sense
and have a
related theoretical development. We give one example here. Recall the
EBEL function (\ref{eqn:EBEL}) for the mean $L_n(\mu)$, $\mu\in
\mathbb{R}^d$,
involves centered block sums $T_{i,\mu} = w(i/n) \sum_{j=1}^{i}(X_j -
\mu)$, $i=1,\ldots,n$, based on
blocks $\{ (X_1), (X_1,X_2),\ldots,\break  (X_1,\ldots,X_n)\}$. Reversed
blocks for example, given by $\{ (X_n), (X_n,X_{n-1}),\ldots,\break
(X_n,\ldots,X_1)\}$, can also be additionally incorporated by defining
further block sums\vadjust{\goodbreak} $T_{n+i,\mu} = w(i/n)\times  \sum_{j=1}^{i}(X_{n-j+1} -
\mu
)$, $i=1,\ldots,n,$ and a corresponding EBEL function
\[
\tilde{L}_{n}(\mu) = \sup\Biggl\{ \prod
_{i=1}^{2n} p_i \dvtx p_i
\geq0, \sum_{i=1}^{2n}p_i = 1,
\sum_{i=1}^{2n} p_i
T_{i,\mu} = 0_d \Biggr\}
\]
and ratio $\tilde{R}_n(\mu)= (2n)^{-2n}\tilde{L}_n(\mu) $. At the true
mean $\mu_0\in\mathbb{R}^d$, the log-ratio
$-\log\tilde{R}_n(\mu_0)= \sum_{i=1}^{2n}\log[1+ \tilde{\lambda
}_{n,\mu_0}^\prime T_{i,\mu_0}]$ can similarly be re-written in terms
of a Lagrange multiplier $\tilde{\lambda}_{n,\mu_0}\in\mathbb{R}^d$
satisfying $0_d=\sum_{i=1 }^{2n} T_{i,\mu_0}/[1+ \tilde{\lambda
}_{n,\mu
_0}^\prime T_{i,\mu_0}] $.
The EL distributional results of the previous subsections
then extend in a natural manner, as described below for the mean
inference case; cf. Theorem \ref{theorem2}. For $0 \leq t \leq1$,
recall $f(t)=w(t)B(t)$ (cf. Theorem \ref{theorem1}), for $B(t) = (B_1(t),\ldots
,B_d(t))^\prime$ denoting a vector of i.i.d. copies $B_1(t),\ldots
,B_d(t)$ of standard Brownian motion on $[0,1]$,
and define additionally $\tilde{f}(t)= w(t)[B(1)-B(1-t)]$.

\begin{theorem} \label{theorem4}
Under assumptions \textup{(A.1)--(A.2)}, as $n \rightarrow\infty$,
\[
n^{1/2} \Sigma^{1/2} \tilde{\lambda}_{n,\mu_0}
\stackrel{d} {\rightarrow} \tilde{Y}_d,\qquad  -\frac{1}{n} \log
\tilde{R}_n(\mu_0) \stackrel{d} {\rightarrow} -
\tilde{g}_d(\tilde{Y}_d)\in(0, \infty)
\]
for a $\mathbb{R}^d$-valued random vector $\tilde{Y}_d$ defined as the
unique minimizer of
\[
\tilde{g}_d(a) \equiv-\int_0^1
\log\bigl( 1 + a^\prime f(t)\bigr)\,dt -\int_0^1
\log\bigl( 1 + a^\prime\tilde{f}(t)\bigr)\,dt
\]
for $a\in\overline{K}_d\equiv\{ y \in\mathbb{R}^d \dvtx\min_{0
\leq t
\leq1} (1 + y^\prime f(t)) \geq0,
\min_{0 \leq t \leq1} (1 + y^\prime\tilde{f}(t)) \geq0\}$.
\end{theorem}
As in Theorem \ref{theorem2}, the limit law of the log EBEL ratio above is similarly
distribution-free and easily simulated from Brownian motion.
The main difference between Theorems \ref{theorem2} and \ref{theorem4}
is that the reversed data blocks in the EL construction
contribute a further integral component based on (reversed) Brownian
motion in the limit. Straightforward analog versions of
Theorem~\ref{theorem1} [regarding $\tilde{Y}_d$ and $\tilde{g}_d(\cdot)$] as well
as Corollary \ref{co1} and Theorem \ref{theorem13} [with respect
to $\tilde{R}_n(\cdot)$] also hold; we state these in the
supplementary materials for completeness.

\section{Numerical studies}\label{sec3}
Here we summarize the results of a simulation study
to investigate the performance of the EBEL method, considering the coverage
accuracy of confidence intervals (CIs) for the process mean.
We considered several real-valued ARMA processes, allowing a variety
of dependence structures with ranges of weak and strong dependence,
defined with respect to an underlying i.i.d. centered $\chi
_1^2$-distributed innovation series; these
processes appear in Table~\ref{tab2} in the following. Other i.i.d. innovation
types (e.g., normal, Bernoulli, Pareto) produced qualitatively similar results.

For each process, we generated 2000 samples of size
$n=250,500, 1000$ for comparing the coverage accuracy of 90\% CIs from
various EL procedures.
We applied the EBEL method with forward\vadjust{\goodbreak} expansive data blocks, as in
Section~\ref{sec2.1}, as
well as forward/backward data blocks, as in Section~\ref{sec2.4}; we denote
these methods as EBEL1/EBEL2, respectively,
in summarizing results. In addition to a constant weight $w(t)=1$, we
implemented these methods
with several other choices of weight functions $w(t)$ on $[0,1]$, each
down-weighting the initial (smaller) data blocks
in the EBEL construction and differing in their shapes. The resulting
coverages were very similar
across nonconstant weight functions and we provide results for two
weight choices: linear $w(t)=t$
and cosine-bell $w(t)=[1-\cos(2\pi t)]/2$.
Additionally, for each weight function $w(t)$, the limiting
distribution of the EBEL ratio was approximated by 50,000
simulations to determine its 90th percentile for calibrating intervals,
as listed in Table~\ref{tab1} with
Monte Carlo error bounds.

%
\begin{table}
\caption{Approximated 90th percentiles of the limit
law of the log-EBEL ratio [$-g_1(Y_1)$ under Theorem \protect\ref{theorem2}
for EBEL1 and $-\tilde{g}_1(\tilde{Y}_1)$ under Theorem \protect\ref{theorem4}
for EBEL2] for weight functions $w(t)$. Approximation $\pm$
parenthetical quantity gives a 95\% CI for
true percentile}\label{tab1}
\begin{tabular*}{\textwidth}{@{\extracolsep{\fill}}lcc@{}}
\hline
$\bolds{w(t)$, $t\in[0,1]}$ & \multicolumn{1}{c}{$\bolds{-g_1(Y_1)}$} &   \multicolumn{1}{c@{}}{$\bolds{-\tilde{g}_1(\tilde{Y}_1)}$}\\
\hline
$w(t)=1$ &2.51 (0.03) &2.50 (0.03)\\
$w(t)=t$ &5.64 (0.09) &4.37 (0.06)\\
$w(t)= (1-\cos(2\pi t))/2$ &7.00 (0.15) & 3.42 (0.09)\\
\hline
\end{tabular*}
\end{table}

For comparison, we also include coverage results for the standard BEL
method with
OL blocks (denoted as BEL). Kitamura \cite{K97} (page 2093) considered
a block order $n^{1/3}$ for BEL
as the method involves a block-based
variance estimator in its asymptotic studentization mechanics (see
Section~\ref{sec4}), which is asymptotically
equivalent to the Bartlett kernel spectral density estimator at zero
having $n^{1/3}$ at its optimal block/lag order; cf. \cite{P03}.
Based this correspondence, we considered two data-driven
block selection rules from the spectral kernel literature, which
estimate the coefficient
$\hat{C}$ in the theoretical optimal block length expression $C n^{1/3}$
known from spectral estimation.
One block estimation approach (denoted FTK) is based on flat-top
kernels and results in block estimates
for BEL due to a procedure in Politis and White \cite{PW04}, page 60;
we used a flat-top kernel bandwidth $n^{1/5}$ for generally consistent
estimation as described in \cite{PW04}.
The second block estimation approach (denoted AAR) is due to
Andrews \cite{A91}, pages 834--835, producing block estimates for BEL
based on bandwidth estimates for the Bartlett spectral kernel assuming
an approximating AR(1)
process.

%
\begin{table}
\caption{Coverage percentages of 90\% intervals for
the process mean over several ARMA processes
(with listed AR/MA components) and sample sizes $n$. EBEL1/EBEL2 use
constant $w(t)=1$ or linear $w(t)=t$ weights; BEL uses
FTK or AAR data-based block selections. [$\operatorname{MA}(1)^*$ has a discrete
component $X_t = \varepsilon_{t} + 0.5\mathbb{I}(\varepsilon_{t-1}<
\chi
^2_{1,0.8})-1.4$, i.i.d. $\varepsilon_t\sim\chi^2_1$]}\label{tab2}
\begin{tabular*}{\textwidth}{@{\extracolsep{\fill}}lcccccccccc@{}}
\hline
& && \multicolumn{2}{c}{\textbf{EBEL1,} $\bolds{w(t)}$} && \multicolumn{2}{c}{\textbf{EBEL2,}
$\bolds{w(t)}$}&& \multicolumn{2}{c@{}}{\textbf{BEL}}\\[-6pt]
& && \multicolumn{2}{c}{\hrulefill} && \multicolumn{2}{c}{\hrulefill}&& \multicolumn{2}{c@{}}{\hrulefill}\\
\textbf{Process} & $\bolds{n}$ && \textbf{1} & $\bolds{t}$&& \textbf{1} &$\bolds{t}$&& \textbf{FTK} &\textbf{AAR}\\
\hline
$\operatorname{MA}(2)$ &\phantom{0}250&&90.6 &91.1 && 91.4 &91.4 &&93.7 &98.3 \\
$0.4, -0.6$&\phantom{0}500&&91.0 &91.2 && 91.7 &91.5 && 93.4 &98.0 \\
&1000&& 90.0& 90.0 && 90.4 &90.0 && 90.6 &96.6\\[3pt]
$\operatorname{MA}(1)^*$&\phantom{0}250&& 87.4 &89.4 && 90.2 &90.5 && 91.3 & 94.2\\
&\phantom{0}500&& 89.4 & 90.8 &&90.4 &90.2 && 90.9 &92.7\\
&1000&&89.6 &89.8 && 90.8 & 90.2 && 91.3 &92.9 \\[3pt]
$\operatorname{MA}(3)$&\phantom{0}250&& 87.4 &88.5 && 90.4& 90.8 && 93.6 &92.7\\
$-1,-1,-1$&\phantom{0}500&& 87.8 &88.6 && 90.0 &90.2 &&93.4 & 92.0 \\
&1000&& 89.7 &89.2 && 89.2 &89.8 && 92.2 &91.9 \\[3pt]
$\operatorname{ARMA}(1,2$)&\phantom{0}250&& 84.4 &86.0 && 89.1& 89.8 && 93.8& 94.6\\
$0.9,-0.6,-0.3$&\phantom{0}500&& 87.2 &88.7 && 90.4& 90.3 && 95.5 &95.2 \\
&1000 && 89.4 &89.9 && 91.6& 91.6 && 95.6& 96.2 \\[3pt]
$\operatorname{AR}(1)$ &\phantom{0}250&&89.2 &90.0&& 92.0& 91.4&& 95.8 &91.8\\
$-0.7$ &\phantom{0}500&& 89.4& 90.6 && 90.9& 90.8 && 95.2& 91.0\\
&1000&& 90.4& 90.2 && 90.4& 90.8 && 92.4 &92.0\\[3pt]
$\operatorname{AR}(1)$ &\phantom{0}250&&67.0 &70.5 && 79.0& 80.0 &&61.1 &76.4 \\
$0.9$ &\phantom{0}500&& 73.4 &77.0 && 82.4& 83.4 &&66.0& 81.4 \\
&1000&& 77.4 &80.1 && 86.2& 87.2 && 74.6& 85.6 \\[3pt]
$\operatorname{ARMA}(1,1)$ &\phantom{0}250&&79.5 &81.8 && 86.3& 86.2 && 81.0 &80.2 \\
$0.7,-0.5$ &\phantom{0}500 && 82.0 &84.6 && 86.3& 86.9&& 82.2 &82.0 \\
&1000 && 85.0& 87.0 && 87.9& 89.0 && 85.4& 84.0 \\[3pt]
$\operatorname{ARMA}(2,2)$ &\phantom{0}250&&78.3 &81.0 && 84.0 &84.6 && 77.2& 73.0 \\
$0.3,0.3,-0.3,-0.1$ &\phantom{0}500 && 81.5 &83.6 && 86.2 &87.2&& 81.0& 74.4 \\
&1000 && 84.4 &85.4 && 88.4 &88.7 && 84.7& 75.3 \\[3pt]
$\operatorname{ARMA}(2,2)$ &\phantom{0}250&&81.2 &83.9 && 85.5& 86.2 && 79.4 &81.8 \\
$0.5,0.3,0.3,-0.9$ &\phantom{0}500 && 84.2& 86.0 && 87.4 &88.0 && 82.8 &84.6 \\
&1000 && 85.4& 86.2 && 88.0& 88.2 && 84.0& 85.5 \\[3pt]
$\operatorname{MA}(2)$ &\phantom{0}250&&83.2 &85.0 && 87.4& 87.6 && 86.0 &79.2 \\
$0.1,2$ &\phantom{0}500 && 84.6 &86.0 && 87.5 &88.6&& 86.8& 81.2 \\
&1000 && 86.2 &87.2 && 89.2 &90.2 && 87.5 &80.4 \\
\hline
\end{tabular*}
\end{table}

Table~\ref{tab2} lists the realized coverage accuracy of 90\% EL CIs for the
mean. From the table, the linear weight
function $w(t)=t$ generally produced slightly more accurate coverages
for both EBEL1/EBEL2 methods than the constant weight $w(t)=1$;
additionally and interestingly, despite their shape differences, the
coverage rates for both the linear and cosine-bell weight functions
closely matched (to the extent that we defer the cosine-bell results to
the supplementary materials~\cite{NBL13}).
For all sample sizes and processes in Table~\ref{tab2}, the EBEL2 method with
linear weight typically and consistently emerged
as having the most accurate coverage properties, often exhibiting less
sensitivity
to the underlying dependence while most closely achieving the nominal
coverage level.
Additionally, linear weight-based EBEL1 generally performed similarly
to, or somewhat better than, the \textit{best} BEL method based on
a data-driven block selection from among the FTK/AAR block rules and,
at times, much better
than the worst performer among the BEL methods with estimated blocks.
Note as well that,
while that the two block selection rules for BEL can produce similar
coverages, their relative effectiveness
often depends crucially on the underlying process, with no resulting
clear best block selection for BEL.
In the case of the strong positive $\operatorname{AR}(1)$ dependence model in Table~\ref{tab2},
the AAR block selection for BEL performed well
(i.e., better
than EBEL1 or BEL/FTK approaches), but similar advantages in coverage
accuracy did not necessarily carry over to other
processes.
In particular, for a process not approximated well by an $\operatorname{AR}(1)$ model,
the BEL coverage rates from AAR block estimates may exhibit extreme
over- or under-coverage
under negative or positive dependence, respectively, and FTK block
selections for BEL may prove better.

%
%
%
%

\begin{figure}

\includegraphics{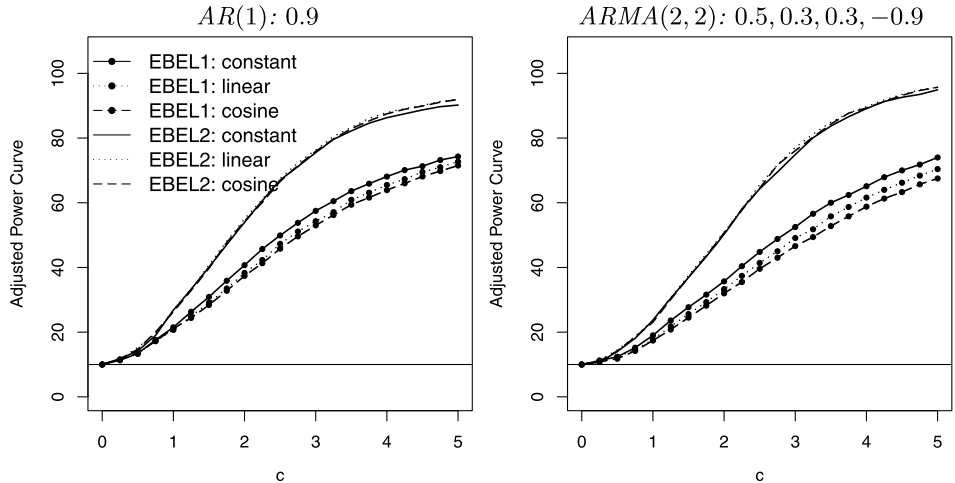}

\caption{Adjusted power curves for tests at 10\% level using
EBEL1/EBEL2 methods with constant, linear and cosine-bell weight
functions (sample size $n=500$).}\label{fig2}
\end{figure}

\begin{figure}

\includegraphics{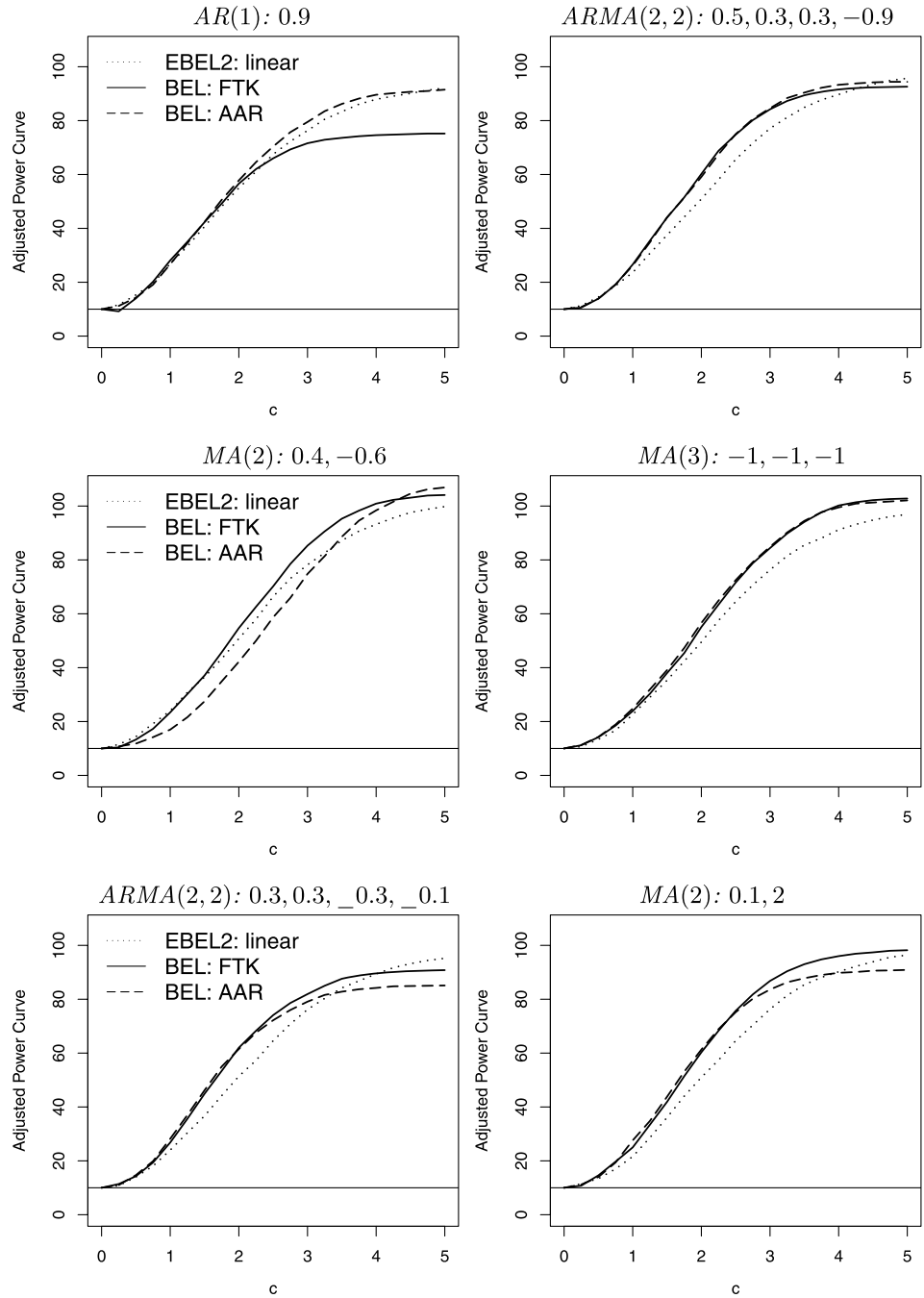}

\caption{Adjusted power curves for tests at 10\% level using BEL with
FTK/AAR block selections and linear weight-based EBEL2 (sample size
$n=500$).}\label{fig3}
\end{figure}

Because of the blocking scheme in EBEL method and some of the method's
other connections to fixed-$b$ asymptotics (see Section~\ref{sec4}), one might
anticipate that
there exist trade-offs in coverage accuracy (i.e., good size control
properties) at the expense of power
in testing, a phenomenon also associated with fixed-$b$ asymptotics; cf.
\cite{BKV01,SPJ08}.
This does seem to be the case. To illustrate, for various sample sizes
$n$ and processes, we approximated power curves
for EBEL/BEL tests at the 10\% level (based on the 90th percentile of
the associated null limit law) along a sequence of
local alternatives $c_n = \mu_0+ n^{-1/2}\Sigma^{1/2}c$,
$c=0,0.25,\ldots,5$ where $\Sigma=\sum_{k=-\infty}^{\infty}\cov
(X_0,X_k)$; for example, with EBEL1, the power curves correspond to the
rejection probabilities $P( -n^{-1}\log R_n(c_n) > q_{0.90} )$ where
$q_{0.90}$ is a percentile from Table~\ref{tab1}.
The alternative sequence $c_n$ was formulated to make power curves
roughly comparable across processes with varying sample sizes, and so
that the power curves can be plotted as a function of $c=0,0.25,\ldots
,5$; for instance, by Corollary \ref{co1}(i), the asymptotic power curve of
EBEL1 will be a function of $c$, as will the curves for BEL/EBEL2.
Figures~\ref{fig2} and \ref{fig3} display \textit{size adjusted} power curves (APCs) for
samples of size $n=500$ based on 2000 simulations (curves are similar
for $n=250,1000$ with additional results given in \cite{NBL13}). If a
percentage $\hat{\alpha}_n$ denotes the actual size of the test
for a given method and process (i.e., $\hat{\alpha}_n= 100\% -{}$coverage percentage in Table~\ref{tab2}), the APC
is
calibrated
to have size 10\% by vertically shifting the true power curve by $10\%
- \hat{\alpha}_n$; this allows the shapes of power curves to be more
easily compared across methods. Figure~\ref{fig2} shows APCs for EBEL1/EBEL2 methods,
where EBEL2 curves exhibit more power apparently as a result
of combining two data block sets (i.e., forward/backward) in the EBEL
construction rather than one; additionally, while EBEL2 power curves
are quite similar across different weights, EBEL1 curves exhibit
slightly more power
for the constant weight function. Figure~\ref{fig3} shows APCs in comparing the
linear weight-based EBEL2 method with
BEL methods based on FTK/AAR block estimates. The APCs for EBEL2
generally tend to be
smaller than those of BEL, though the APC of a block estimate-based
BEL may not always
dominate the associated curve of EBEL2.

%

%
%
%

%

%

\section{Conclusions}\label{sec4}
The proposed expansive block empirical likelihood\break  (EBEL) is a type of
variation on standard blockwise empirical likelihood (BEL) for time
series which, instead of using a fixed
block length $b$ for a given sample size $n$, involves a nonstandard blocking
scheme to capture the dependence structure. While the coverage accuracy
of standard BEL methods can depend intricately
on the block choice $b$ (where the best $b$ can vary with the
underlying process),
the EBEL method does not involve this type of block selection. As
mentioned in the \hyperref[sec1]{Introduction}, we also anticipate
that the EBEL method
will generally have better rates
of coverage accuracy compared to BEL. The simulations of
Section~\ref{sec3} lend support to this notion, along with suggesting that the EBEL
can be less sensitive to the strength of the underlying time dependence.
While asymptotic coverage
rates for BEL methods remain to be determined, we may offer the
following heuristic based on analogs drawn to so-called ``fixed-$b$ asymptotics''
(cf. Keifer, Vogelsang and Bunzel \cite{KVB00}; Bunzel et al. \cite
{BKV01}; Kiefer and Vogelsang \cite{KV02}),
or related ``self-normalization'' (cf. Lobato \cite{L01}; Shao \cite
{S10}) schemes.

In asymptotic expansions of log-likelihood statistics from standard BEL
formulations, the data blocks serve
to provide a type of block-based variance estimator (cf. \cite{C86,PR93})
for purposes of normalizing scale and obtaining chi-square limits
for log-BEL ratio statistics. Such variance estimators are consistent,
requiring block sizes $b$ which grow at a smaller rate than the sample
size $n$ (i.e., $b^{-1}+b/n\rightarrow\infty$ as $n\rightarrow\infty
$) and are known to
have equivalences to variance estimators formulated as lag window
estimates involving kernel functions and bandwidths $b$ with similar
behavior to block lengths $b^{-1}+b/n\rightarrow\infty$ (cf. \cite
{K89,P03}). That is,
standard BEL intervals have parallels with normal theory intervals
based on normalization with consistent lag window estimates.
However, considering hypothesis testing with sample means, for example,
there is some numerical and theoretical evidence (cf. \cite
{BKV01,SPJ08}) that normalizing scale with
inconsistent lag window estimates having fixed bandwidth ratios (e.g.,
$b/n =C$ for some $C\in(0,1]$) results in
better size and lower power
compared to normalization with consistent ones, though the former case requires
calibrating intervals with nonnormal limit laws.
Shao \cite{S10}, Section~2.1, provides a nice summary of these points
as well as the form of some of these distribution-free limit laws,
which typically
involve ratios of random variables defined by Brownian motion; cf.
\cite{KV02}.
While the EBEL method is not immediately analogous to normalizing with
inconsistent variance estimators (as mentioned in Section~\ref{sec2.1}, the
usual EL expansions do not hold
for EBEL), there are parallels in that the EBEL
method does not use block lengths satisfying standard bandwidth
conditions (cf. Section~\ref{sec2.1}), its blocking scheme itself
appears in self-normalization literature (cf. Shao \cite{S10}, Section~2)
and confidence region calibration involves nonnormal limits based on
Brownian motion.
This heuristic in the mean case suggests that better coverage rates
(and lower power) associated
with fixed-$b$ asymptotics over standard normal theory asymptotics
may be anticipated to carry over to comparisons of EBEL to standard
BEL formulations.

\section{Proofs of main results}\label{sec5}

To establish Theorem \ref{theorem1}, we first require a lemma regarding a standard
Brownian motion.
For concreteness, suppose $B(t)\equiv B(\omega,t) = (B_1(\omega
,t),\ldots,B_d(\omega,t))^\prime$, $\omega\in\Omega$, $t\in[0,1]$
is a random $\mathcal{C}_d[0,1]$-valued element defined on some
probability space $(\Omega,\mathcal{F},P)$,
where $B_1,\ldots,B_d$ are again distributed as i.i.d. copies of
standard Brownian motion on $[0,1]$.
In the following, we use the basic fact that each $B_i(\cdot)$ is
continuous on $[0,1]$ with probability~1 (w.p.1)
along with the fact that increments of standard Brownian motion are
independent; cf. \cite{F83}.

%
\begin{Lemma}\label{le1} With probability 1, it holds that:
\begin{longlist}[(iii)]
\item[(i)] $ {\min_{0 \leq t < \epsilon} a^\prime B(t) < 0 < \max_{0
\leq t < \epsilon} a^\prime B(t)}$ for all $\epsilon>0$ and $a\in
\mathbb{R}^d$,\break  $\|a\|=1$.
\item[(ii)] $0_d$ is in the interior of the convex hull of $B(t)$, $0
\leq t \leq1$.
\item[(iii)] There exists a positive random variable $M$ such that, for
all $a \in\mathbb{R}^d$, it holds that $ {\min_{0 \leq t \leq1}
a^\prime B(t) \leq-M \|a\|}$ and $ {M\|a\| \leq\max_{0 \leq t \leq1}
a^\prime B(t)}$.
\item[(iv)] If assumption \textup{(A.1)} holds in addition, \textup{(i)}, \textup{(ii)}, \textup{(iii)} above
hold upon replacing $B(t)$ with $f(t)=w(t)B(t)$, $t\in[0,1]$.
\end{longlist}
\end{Lemma}
\begin{pf}
For real-valued Brownian motion,
it is known that\break  $\min_{0 \leq t < \epsilon} B_i(t) <  0 < \max_{0
\leq t < \epsilon} B_i(t)$ holds for all $\epsilon>0$ w.p.1. (cf.
\cite
{F83}, Lemma 55); we modify the proof of this. Let $\{t_n\}\subset
(0,1)$ be a decreasing sequence where $t_n \downarrow0$ as
$n\rightarrow\infty$. Pick and fix $c_1,\ldots,c_d\in\{-1,1\}$, and
define the event $A_n \equiv A_{n,c_1,\ldots,c_d} =\{\omega\in\Omega
\dvtx
c_i B_i(\omega,t_n) >0, i=1,\ldots,d\}$. Then $P(A_n) = 2^{-d}$
for all $n\geq1$ by normality and independence. As the events $B_n =
\bigcup_{k=n}^\infty A_k$, $n \geq1$, are decreasing, it holds that
\[
P \Biggl( \bigcap_{n=1}^\infty
B_n \Biggr) = \lim_{n\to\infty}P(B_n) \geq\lim
_{n\to\infty}P(A_n) = 2^{-d}.
\]
Since $\bigcap_{n=1}^\infty B_n$ is a tail event generated by the
independent random variables $B_i(t_1)- B_i(t_2),
B_i(t_2)-B_i(t_3),\ldots$ for $i=1,\ldots,d$ [i.e., increments of
Brownian motion are independent and $B_i(0)=0$), it follows from
Kolmogorov's 0--1 law that $1 =P ( \bigcap_{n=1}^\infty B_n )
= P ( \mbox{$A_n$ infinitely often (i.o.)} ]$. Hence
$P ( \mbox{$A_{n,c_1,\ldots,c_d}$}\break \mbox{i.o. for any $c_i\in\{1,-1\}$,
$i=1,\ldots,d$} ) =1$ must hold, which implies part (i).

For part (ii), if $0_d$ is not in the interior convex hull of $B(t)$,
$t\in[0,1]$, then the supporting/separating hyperplane theorem would
imply that, for some $a\in\mathbb{R}^d$, $\|a\|=1$, it holds that
$a^\prime B(t)\geq0$ for all $t\in[0,1]$, which contradicts part (i).

To show part (iii), we use the events developed in part (i) and define
$n_{c_1,\ldots,c_d} = \min\{n \dvtx\mbox{$A_{n,c_1,\ldots,c_d}$
holds}\}
$. Define $M = \min\{ |B_i(t_{n_{c_1,\ldots,c_d}})| \dvtx c_1,\ldots,c_d
\in\{-1,1\},  i=1,\ldots,d\}>0$. For $a=(a_1,\ldots,a_d)^\prime\in
\mathbb{R}^d$, let $c_i^{a} = \max\{-\operatorname{sign}(a_i),1\}$,
$i=1,\ldots,d$. Then $a^\prime B(t_{n_{c_1^a,\ldots,c_d^a}}) = - \sum
_{i=1}^d |a_i B_i(t_{n_{c_1^a,\ldots,c_d^a}})| \leq-M \|a\| $, and
likewise $a^\prime B(t_{n_{-c_1^a,\ldots,-c_d^a}}) = \sum_{i=1}^d |a_i
B_i(t_{n_{-c_1^a,\ldots,-c_d^a}})| \geq M \|a\|$. This estab-\break lishes~(iii).

Part (iv) follows from the fact that $w(t)>0$ for $t \in(0,c)$, and we
may take the positive sequence $\{t_n\}\subset(0,c)$ in the proof of
part (i). Then the results for $B(t)$ imply the same hold upon
substituting $f(t)=w(t)B(t)$, $t\in[0,1]$.
\end{pf}

\begin{pf*}{Proof of Theorem \ref{theorem1}} The set $K_d = \{a\in\mathbb
{R}^d \dvtx\min_{0 \leq t \leq1} (1+a^\prime f(t))>0\}$
is open, bounded and convex (w.p.1), where boundedness follows from
Lemma \ref{le1}(iii), (iv). Likewise, the closure $\overline{K}_d = \{a\in
\mathbb
{R}^d \dvtx\min_{0 \leq t \leq1} (1+a^\prime f(t))\geq0\}$ is
convex and
bounded. Since $\min_{0 \leq t \leq1} (1+a^\prime f(t))$ is a
continuous function in $a\in\mathbb{R}^d$, one may apply the dominated
convergence theorem (DCT)
[with the fact that $\min_{0 \leq t \leq1}(1+a^\prime f(t)]$ is
bounded away from 0 on closed balls inside $K_d$ around $a$)
to show that partial derivatives of $g_d(\cdot)$ at $a\in K_d$ (of all
orders) exist, with first and second partial derivatives given by
\[
\frac{\partial g_d(a)}{\partial a} = - \int_0^1
\frac{f(t)}{ 1+a^\prime
f(t)}\,dt, \qquad\frac{\partial^2 g_d(a)}{\partial a\, \partial a^\prime} =
\int_0^1
\frac{f(t)f(t)^\prime}{ [1+a^\prime f(t)]^2}\,dt.
\]
Because $\int_0^1 f(t)f(t)^\prime \,dt$ is positive definite by
Lemma \ref{le1}(i), (iv) and the continuity of $f$,
the matrix $\partial^2 g_d(a)/\partial a \,\partial a^\prime$ is also
positive definitive for all $a\in K_d$, implying
$g_d$ is strictly convex on $K_d$. By Jensen's inequality, it also
holds that $g_d$ is convex on $\overline{K}_d$.

Note for $a\in\overline{K}_d$, $g_d(a) \geq- \int_0^1 \log(1 +
\sup_{a\in\overline{K}_d}\|a\|\cdot\sup_{0 \leq t \leq1} \|f(t)\|
) >
-\infty$ holds, so that $I \equiv\inf_{a\in\overline{K}_d}g_d(a)$
exists. Additionally, $0_d \in K_d$ with $g_d(0_d)=0$ and
$\partial g_d(0_d)/\partial a = - \int_0^1 f(t)\,dt $, where the
components of $\int_0^1 f(t)\,dt$ are all nonzero (w.p.1) by normality
and independence; by the continuity
of partial derivatives on the open set $K_d$, there then exists $\bar
{a} \in K_d$ such that $\bar{a}^\prime\int_0^1 f(t)\,dt>0$ holds with
the components of $- \int_0^1 f(t)$ and $\partial g_d(\bar
{a})/\partial a$ having the same sign. By strict convexity,
$g_d(0_d)-g_d(\bar{a})> [\partial g_d(\bar{a})/\partial a]^\prime
(0_d-\bar{a})>0$ follows, implying
$I< 0$ and $I = \inf_{a\in\tilde{K}_d}g_d(a)$
for the level set $\tilde{K}_d \equiv\{a \in\overline{K}_d\dvtx
g_d(a)\leq0\}$.

Then, there exists a sequence $a_n \in\tilde{K}_d $ such that
$g_d(a_n)< I + n^{-1}$ for $n \geq1$. Since $\{a_n\}$ is bounded, we
may extract a subsequence such that $a_{n_k}\rightarrow Y_d \in\tilde
{K}_d$, for some $Y_d \in\tilde{K}_d$. Pick $\delta\in(0,1)$. Then, by
the DCT,
\begin{eqnarray*}
\underline{\lim}\, g_d(a_{n_k}) &\geq& \underline{\lim} \int
_{\{t:
a_{n_k}^\prime f(t)> -1 + \delta\}} -\log\bigl(1 + a_{n_k}^\prime f(t)
\bigr)\,dt
\\
& = & \int_{\{t: Y_d^\prime f(t)> -1 + \delta\}} -\log\bigl(1
+Y_d^\prime
f(t) \bigr)\,dt
\\
&=& g_d(Y_d) + \int_{\{t: Y_d^\prime f(t) \leq-1 + \delta\}} \log
\bigl(1
+Y_d^\prime f(t) \bigr)\,dt.
\end{eqnarray*}
Note that because $g_d(Y_d)\in(-\infty,0]$, it follows that $-\int
_{\{
t: Y_d^\prime f(t) <0\}} \log(1 +Y_d^\prime f(t) )\,dt<\infty$ and $\{
t\in[0,1]\dvtx Y_d^\prime f(t)=-1\}$ has Lebesgue measure zero. Hence, the
DCT yields
\[
\lim_{\delta\rightarrow0}-\int_{\{t: Y_d^\prime f(t) \leq-1 +
\delta
\}} \log\bigl(1
+Y_d^\prime f(t) \bigr)\,dt=0.
\]
Consequently,
\[
I \geq\overline{\lim}\, g_d(a_{n_k})\geq\underline{\lim}\,
g_d(a_{n_k}) \geq g_d(Y_d) \geq I,
\]
establishing the existence of a minimizer $Y_d$ of $g_d$ on $\overline
{K}_d$ such that $-\infty< I=g_d(Y_d)<0$.

For part (ii) of Theorem \ref{theorem1}, note $y_n = (1-n^{-1}) Y_d + n^{-1}0_d \in
K_d$, $n\geq1$, by convex geometry, as $K_d$ is the convex interior of
$\overline{K}_d$. Then $g_d(y_n) \leq(1-n^{-1}) g_d(Y_d)$ holds by
convexity of $g_d$ and $g_d(0_d)=0$, implying $0 \leq n[g_d(y_n) -
g_d(Y_d)] \leq-g_d(Y_d)<\infty$, from which it follows that
$g_d(y_n)\rightarrow g_d(Y_d)$
and, by the mean value theorem,
\[
0 \leq n\bigl[g_d(y_n) - g_d(Y_d)
\bigr]=\int_0^1 \frac{Y_d^\prime f(t)}{1 + c_n
Y_d^\prime f(t)} \,dt
\leq-g_d(Y_d)
\]
holds
for some $(1-n^{-1})<c_n < 1$ [note $c_n Y_d \in K_d$ so $\min_{0 \leq
t \leq1} (1 + c_n Y_d^\prime f(t))>0$ for all $n$]; the latter implies
$0 \leq\int_{\{t: Y_d^\prime f(t)<0\}} -Y_d^\prime f(t)/[1 + c_n
Y_d^\prime f(t)] \,dt \leq\int_{\{t: Y_d^\prime f(t)>0\}} Y_d^\prime
f(t)<\infty$ so that Fatou's lemma yields
\[
0 \leq\int_{\{t: Y_d^\prime f(t)<0\}} - \frac{Y_d^\prime f(t)}{1 +
Y_d^\prime f(t)} \,dt <\infty
\]
as $n\rightarrow\infty$, and consequently $\int_0^1 1/[1 +
Y_d^\prime
f(t)]\,dt <\infty$. We may then apply the DCT to find
\[
\lim_{n\to\infty} \int_0^1
\frac{Y_d^\prime f(t)}{1 + c_n Y_d^\prime
f(t)} \,dt = \int_0^1
\frac{Y_d^\prime f(t)}{1 + Y_d^\prime f(t)} \,dt \in[0,\infty).
\]
Also by convexity and $0_d\in K_d$, $0>g_d(Y_d) - g_d(0_d) > [\partial
g_d(0_d)/\partial a]^\prime(Y_d-0_d)$ holds (w.p.1), implying
$Y_d^\prime\int_0^1 f(t)\,dt >0$ from $\partial g_d(0_d)/\partial a=-
\int_0^1 f(t)\,dt$. This establishes part (ii) of Theorem \ref{theorem1}.

To show uniqueness of the minimizer, we shall construct sequences with
the same properties
in the proof of part (ii) above. Suppose $x \in\tilde{K}_d$ such that
$g_d(x)=I=g_d(Y_d)$. Defining $x_n =(1-n^{-1})x + n^{-1}0_d \in K_d$
and $y_n = (1-n^{-1}) Y_d + n^{-1}0_d \in K_d$ for $n\geq1$, by
convexity we have $0 \geq g_d(x)-g_d(y_n)> [\partial g_d(y_n)/\partial
a]^\prime(x-y_n)$, so that taking limits yields $0 \geq- \int_0^1
(x-Y_d)^\prime f(t)/[1 + Y_d^\prime f(t)] \,dt$, and, by symmetry,
$0 \geq- \int_0^1 (Y_d-x)^\prime f(t)/[1 + x^\prime f(t)] \,dt$ as well.
Adding these terms gives
\[
0 \geq\int_0^1 \frac{[(x-Y_d)^\prime f(t)]^2 }{(1 + x^\prime f(t))(1 +
Y_d^\prime f(t))} \,dt,
\]
implying that $x=Y_d$ by Lemma \ref{le1}(iv) and
the continuity of $f$.

Finally, to establish part (iii), if $Y_d \in K_d$, then $0_d =
\partial g_d(Y_d)/\partial a = - \int_0^1 f(t)/\break  [1 + Y_d^\prime f(t)]\,dt$
must hold. If there exists another $b \in\overline{K}_d$ satisfying
$\int_0^1 f(t)/[1 + b^\prime f(t)]\,dt =0_d$, then
adding $\partial g_d(Y_d)/\partial a$ to this integral and multiplying
by $(Y_d-b)^\prime$ yields
$0=\int_0^1 [(b-Y_d)^\prime f(t)]^2 /[(1 + b^\prime f(t))(1 +
Y_d^\prime f(t))] \,dt$, implying that $b=Y_d$. Also,
if $ 0_d=\int_0^1 f(t)/[1 + b^\prime f(t)]\,dt = - \partial
g_d(b)/\partial a $ holds for some $b \in K_d$, then strict convexity implies
$g_d(a)-g_d(b) > [\partial g_d(b)/\partial a]^\prime(a-b)=0$ for all $a
\in\overline{K}_d$, implying $b=Y_d$ is the unique
minimizer of $g_d$.
\end{pf*}

\begin{pf*}{Proof of Theorem \ref{theorem2}} Under assumption (A.2), we use
Skorohod's embedding theorem (cf. \cite{VW96}, Theorem 1.1.04) to embed
$\{S_n(\cdot)\}$ and $\{B(\cdot)\}$ in a larger probability space
$(\Omega,\mathcal{F},P)$ such that
$\sup_{0 \leq t \leq1}\| \Sigma^{-1/2} S_n(t)/n^{1/2} -B(t) \|
\rightarrow0$ w.p.1 $(P)$. Defining
$T_n(t)=w(t)S_n(t)$ and $f(t)=w(t)B(t)$, $t\in[0,1]$, the continuity
of $w$ under assumption (A.1) then implies
%
%
\begin{equation}
\label{eqn:Tn} \sup_{0 \leq t \leq1}\biggl\llVert\frac{\Sigma
^{-1/2}T_n(t)}{n^{1/2}}
-f(t) \biggr\rrVert\rightarrow0 \qquad\mbox{w.p.1.}
\end{equation}
Note that $T_{i,\mu_0}= w(i/n)\sum_{j=1}^i(X_j - \mu_0) = T_n(i/n)$,
$i=1,\ldots,n$. By (\ref{eqn:Tn}) and Lemma \ref{le1},
$0_d$ is in the interior convex hull of $\{T_{i,\mu_0}\dvtx i=1,\ldots
,n\}$
eventually (w.p.1) so that $L_n(\mu_0)>0$ eventually (w.p.1). That is,
by Lemma \ref{le1}(iv), there
exists $A\in\mathcal{F}$ with $P(A)=1$ and, for $\omega\in A$, $\min
_{0 \leq t \leq} a^\prime f(\omega,t) \leq-M(\omega)$
and $\max_{0 \leq t \leq} a^\prime f(\omega,t) \geq M(\omega)$ hold
for some $M(\omega)>0$ and all $a\in\mathbb{R}^d$, $\|a\|=1$. Then,
(\ref{eqn:Tn}) implies\vspace*{1pt} $\min_{1\leq i \leq n} a^\prime\Sigma
^{-1/2}T_n(\omega, i/n) < 0 < a^\prime\max_{1\leq i \leq n} \Sigma
^{-1/2}T_n(\omega,\break   i/n)$ holds for all $a\in\mathbb{R}^d$, $\|a\|=1$
eventually, implying $0_d$ is in the interior convex hull of $\{\Sigma
^{-1/2}T_n(i/n)\dvtx i=1,\ldots,n\}$. Hence, eventually (w.p.1) as in
(\ref
{eqn:t0}), we can write
\[
\frac{1}{n} R_n(\mu_0) = -\frac{1}{n}\sum
_{i=1}^n \log\bigl(1 + \lambda
_{n,\mu_0}^{\prime} T_{i,\mu_0}\bigr) = \frac{1}{n}\sum
_{i=1}^n \log\bigl(1 +
\ell_n^{\prime} T_{i,n}\bigr),
\]
where $ T_{i,n} \equiv\Sigma^{-1/2}T_n(i/n)/n^{1/2}$, $i=1,\ldots,n$
and $\ell_n = n^{1/2} \Sigma^{1/2} \lambda_{n,\mu_0}$ and
%
%
\begin{eqnarray}
\label{eqn:Tn2} \min_{i=1,\ldots,n}\bigl(1 + \ell_n^{\prime}
T_{i,n}\bigr)&>&0,\qquad \sum_{i=1}^n
\frac
{1}{n(1 + \ell_n^{\prime}T_{i,n})} =1,
\nonumber
\\[-8pt]
\\[-8pt]
\nonumber
 \sum_{i=1}^n
\frac{T_{i,n}}{n(1
+ \ell_n^{\prime} T_{i,n})}& =&0_d.
\end{eqnarray}

From here, all considered convergence will be pointwise along some
fixed $\omega\in A$ where $P(A)=1$,
and we suppress the dependence of terms $f$, $T_n$, etc. on $\omega$.
Then, (\ref{eqn:Tn2}) [i.e.,
$\min_{i=0,\ldots,n}\ell_n^\prime(\Sigma^{-1/2}T_n(i/n)/n^{1/2}) >-1$]
with (\ref{eqn:Tn}) and Lemma~\ref{le1}(iv) implies that
$\|\ell_n\|$ is bounded eventually. For any subsequence $\{n_j\}$ of
$\{
n\}$, we may extract a further subsequence $\{n_k\} \subset\{n_j\}$
such that $\ell_{n_k}\rightarrow b$ for some $b \in\overline{K}_d$.
For simplicity, write $n_k \equiv k$ in the following. We will show
below that $k^{-1} \log R_k(\mu_0) \rightarrow g_{d}(Y_d)$ and that
$\ell_k \rightarrow Y_d$, where
$Y_d\in\overline{K}_d$ denotes the minimizer of $g_{d}(a)= - \int_0^1
\log(1 + a^\prime f(t))\,dt$, $a\in\overline{K}_d$ under Theorem \ref{theorem1}.
Since the subsequence $\{n_j\}$ is arbitrary, we then have
$n^{-1} \log R_n(\mu_0) \rightarrow g_{d}(Y_d)$ and $\ell_n
\rightarrow Y_d$ w.p.1, implying the distributional convergence in
Theorem \ref{theorem2}.

Define
$Y_\epsilon= (1 - \epsilon)Y_d + \epsilon0_d \in K_d $ (since $0_d
\in K_d$, the interior of $\overline{K}_d$) for $\epsilon\in(0,1)$.
From $Y_\epsilon\in K_d$, $\min_{0 \leq t \leq1} (1 + Y_\epsilon
^\prime f(t))>\delta$ holds for some $\delta>0$ (dependent on
$\epsilon
$) so that $\min_{1 \leq i \leq k} (1 + Y_\epsilon^\prime
T_{i,k})>\delta$ holds eventually by (\ref{eqn:Tn}). Then because
\[
g_{d,k}(a) \equiv- \frac{1}{k} \sum
_{i=1}^k \log\bigl(1+ a^\prime
T_{i,k}\bigr)
\]
is strictly convex on $a\in\{ y\in\mathbb{R}^d \dvtx\min_{1\leq i
\leq
k} (1+y^\prime T_{i,n})>0\}$ with a unique minimizer at $\ell_k$ by
(\ref{eqn:Tn2}) [i.e., $\partial g_{d,k}(\ell_k)/\partial a = 0_d$
holds and strict convexity follows when $k^{-1}\sum_{i=1}^k T_{i,k}
T_{i,k}^\prime$ is positive definite, which holds eventually from
$k^{-1}\sum_{i=1}^k T_{i,k} T_{i,k}^\prime\rightarrow\int_0^1
f(t)f(t)^\prime \,dt$ by (\ref{eqn:Tn}) and the DCT,\vadjust{\goodbreak} with the latter
matrix being positive definite w.p.1 by Lemma \ref{le1}(iv) and continuity of
$f$], we have that
\[
g_{d,k}(Y_\epsilon) \geq g_{d,k}(\ell_k)
= \frac{1}{k} \log R_k(\mu_0).
\]
Define $\bar{g}_{d,k}(a) \equiv- k^{-1}\sum_{i=1}^k \log( 1 +
a^\prime f(i/k))$, $a \in K_d$. Then, by Taylor expansion [recalling
$\min_{0 \leq t \leq1} (1 + Y_\epsilon^\prime f(t))>\delta$, $\min_{1
\leq i \leq k} (1 + Y_\epsilon^\prime T_{i,k})>\delta$],
\begin{eqnarray*}
&&\bigl|g_{d,k}(Y_\epsilon) - \bar{g}_{d,k}(Y_\epsilon)\bigr|
\\
&& \qquad\leq\frac{1}{k}\sum_{i=1}^k
\bigl|Y_\epsilon^\prime\bigl(T_{i,k} -f(i/k)\bigr)\bigr| \biggl(
\frac{1}{1 + Y_\epsilon^\prime T_{i,k}} + \frac{1}{1 + Y_\epsilon
^\prime f(i/k)} \biggr)
\\
&&\qquad\leq \|Y_d\| 2 \delta^{-1} \max_{1 \leq i \leq k}
\bigl\|T_{i,k} -f(i/k) \bigr\| \rightarrow0
\end{eqnarray*}
from (\ref{eqn:Tn}) and Theorem \ref{theorem1}. Also, by the DCT, $\bar
{g}_{d,k}(Y_\epsilon) \rightarrow g_{d}(Y_\epsilon)$ as $k\rightarrow
\infty$. Hence, $g_{d}(Y_\epsilon) \geq\overline{\lim}\,
g_{d,k}(\ell
_k)$ holds and, since $g_{d}( Y_\epsilon) \leq(1-\epsilon) g_{d}(Y_d)$
by convexity and $g_{d}(0_d)=0$, we have, letting $\epsilon\rightarrow
0$, that
%
%
\begin{equation}
\label{eqn:g1} g_{d}(Y_d) \geq\overline{\lim}\,
g_{d,k}(\ell_k).
\end{equation}

Recalling $\ell_k \rightarrow b \in\overline{K}_d$, define
$b_\epsilon
= (1-\epsilon)b + \epsilon0_d \in K_d$, so that
$\min_{0 \leq t \leq1} (1+b_\epsilon^\prime f(t))>0$. Then, $\bar
{g}_{d,k}(b_\epsilon) \rightarrow g_{d}(b_\epsilon)$
by (\ref{eqn:Tn}) and the DCT. And, by Taylor expansion and using
(\ref
{eqn:Tn2}),
\begin{eqnarray*}
&&\overline{\lim} \bigl|g_{d,k}(\ell_k) -
\bar{g}_{d,k}(b_\epsilon)\bigr|
\\
&&\qquad\leq \overline{\lim} \max_{1 \leq i \leq k}\bigl| \ell_k^\prime
T_{i,k} - b_\epsilon^\prime f(i/k)\bigr| \Biggl( 1 +
\frac{1}{k}\sum_{i=1}^k
\frac{1}{1
+ b_\epsilon^\prime f(i/k)} \Biggr)
\\
&&\qquad\leq \epsilon\sup_{0 \leq t \leq1} \bigl|b^\prime f(t)\bigr| \biggl( 1 +
\int_0^1 \frac{1}{1 + b_\epsilon^\prime f(t)}\,dt \biggr)\equiv C(
\epsilon),
\end{eqnarray*}
following from (\ref{eqn:Tn}) and the DCT. Hence we have
%
%
\begin{equation}
\label{eqn:g2} \underline{\lim}\,g_{d,k}(\ell_k) \geq
g_{d}(b_\epsilon) - C(\epsilon).
\end{equation}
We will show below that
%
%
\begin{equation}
\label{eqn:b} \int_0^1 \frac{1}{1 + b^\prime f(t)} \,dt
<\infty
\end{equation}
holds, in which case, $\lim_{\epsilon\rightarrow0}\int_0^1 [1 +
b_\epsilon^\prime f(t)]^{-1}\,dt = \int_0^1 [1 + b^\prime f(t)]^{-1} \,dt
<\infty$ by the DCT and so that $C(\epsilon)\rightarrow0 $ as
$\epsilon\rightarrow0$ [noting
$\sup_{0 \leq t \leq1} |b^\prime f(t)|<\infty$ since $f$ is continuous
and $\overline{K}_d$ is bounded by Theorem \ref{theorem1}].
By Fatou's lemma and the DCT, $\underline{\lim}_{\epsilon\rightarrow0}
g_{d}(b_\epsilon)\geq g_{d}(b)$ holds also. Hence, by (\ref
{eqn:g1})--(\ref{eqn:g2}), we then have
\[
g_{d}(Y_d) \geq\overline{\lim}\, g_{d,k}(
\ell_k) \geq\underline{\lim}\,g_{d,k}(\ell_k)
\geq g_{d}(b) \geq g_{d}(Y_d),
\]
implying $b=Y_d$ by the uniqueness of the minimizer and $\lim_{k\to
\infty} k^{-1}\times \log R_{k}(\mu_0) = g_{d}(Y_d)$.\vadjust{\goodbreak}

To finally show (\ref{eqn:b}), let $A = \{t\in[0,1]\dvtx1 + b^\prime f(t)
\leq d\}$ for some $0<d\leq1/2$ chosen
so that $\{t\in[0,1]\dvtx1 + b^\prime f(t)=d\}$ has Lebesgue measure
zero (since $f$ is continuous). Let $A^c = [0,1]\setminus A$. Using the
indicator function $\mathbb{I}(\cdot)$, define a simple function
\[
h_k(t) \equiv\sum_{i=1}^k
\frac{\ell_k^\prime T_{i,k}}{1 + \ell
_k^\prime T_{i,k}} \mathbb{I} \biggl(t \in\biggl(\frac{i-1}{k}, \frac
{i}{k}
\biggr] \biggr),\qquad  t\in[0,1].
\]
From (\ref{eqn:Tn2}), note that
\[
\int_A h_k(t)\,dt + \int_{A^c}
h_k(t) \,dt = \frac{1}{k} \sum_{i=1}^k
\frac
{\ell_k^\prime T_{i,k}}{1 + \ell_k^\prime T_{i,k}} =0_d.
\]
From (\ref{eqn:Tn}), $\mathbb{I}(t\in A^c) h_k(t) \rightarrow\mathbb
{I}(t\in A^c) b^\prime f(t)/(1+b^\prime f(t))$ [almost everywhere
(a.e.) Lebesgue measure] and
for large $k$, $\mathbb{I}(t\in A^c) |h_k(t)|\leq2C/d$ holds for
$t\in
[0,1]$, since eventually $\max_{1 \leq i \leq k} |\ell_k^\prime
T_{i,k}|$ is bounded by a constant $C>0$ and also $1 + b^\prime f(t) +
(\ell_k^\prime T_{i,k}- b^\prime f(t))>d/2$ for $t\in A^c$, $(i-1)/k< t
\leq i/k$. Then, by the DCT, $\int_{A^c} h_k(t) \,dt \rightarrow\int
_{A^c} b^\prime f(t)/(1+b^\prime f(t)) \,dt$, and for $\delta\in(0,1)$,
note $-\mathbb{I}(t\in A)h_k(t) \geq h_{1,k}(t)$ for
\[
h_{1,k}(t) \equiv\sum_{i=1}^k
\frac{-\ell_k^\prime T_{i,k}}{1 + \ell
_k^\prime T_{i,k} + \delta\mathbb{I}(\operatorname{sign}(\ell_k^\prime
T_{i,k}) <0)} \mathbb{I} \biggl(t \in\biggl(\frac{i-1}{k}, \frac{i}{k}
\biggr]\cap A \biggr).
\]
Since $|h_{1,k}(t)| \leq C/\delta$ and $h_{1,k}(t)\rightarrow-\mathbb
{I}(t\in A) b^\prime f(t) /(1+b^\prime f(t)+\delta) $ (a.e. Lebesgue
measure), by the DCT
\begin{eqnarray*}
0&\leq&\int_A \frac{-b^\prime f(t)}{1+b^\prime f(t) + \delta}\,dt =
\lim
_{k\to\infty} \int_A h_{1,k}(t)\,dt
\leq\lim_{k\to\infty} \int_A
-h_{k}(t)\,dt \\
&=& \int_{A^c} \frac{b^\prime f(t)}{1+b^\prime f(t)}\,dt
\end{eqnarray*}
using $\int_A -h_{k}(t)\,dt = \int_{A^c} h_{k}(t)\,dt$. Letting $\delta
\rightarrow0$, Fatou's lemma gives
\[
0\leq\int_A \frac{-b^\prime f(t)}{1+b^\prime f(t) }\,dt \leq\int
_{A^c} \frac{b^\prime f(t)}{1+b^\prime f(t)}\,dt<\infty.
\]
Because $-b^\prime f(t) \geq1/2$ on $A$, $\int_A [1+b^\prime f(t)
]^{-1}\,dt<\infty$ holds, implying (\ref{eqn:b}).
\end{pf*}

\section*{Acknowledgements}
The authors are very grateful to Editor Runze Li, two Associate Editors
and two referees for constructive and insightful comments which
improved the manuscript, especially the numerical studies.

\begin{supplement}[id=suppA]
\stitle{Additional proofs and results for a nonstandard empirical
likelihood for time series.}
\slink[doi]{10.1214/13-AOS1174SUPP} 
\sdatatype{.pdf}
\sfilename{aos1174\_supp.pdf}
\sdescription{A supplement\vadjust{\goodbreak} \cite{NBL13} provides proofs of the
remaining main results omitted here, namely Corollary~\ref{co1} (properties of
confidence regions), Theorem~\ref{theorem13} (smooth function model results) and
Theorem~\ref{theorem4} (forward/backward block EL version); additional numerical
summaries are included as well.}
\end{supplement}

%


\printaddresses

\end{document}